\documentclass{amsart}
\usepackage{amssymb}
\usepackage{lattice}

\theoremstyle{plain}
\newtheorem{theorem}{Theorem}
\newtheorem{lemma}{Lemma}[section]
\newtheorem{proposition}[lemma]{Proposition}
\newtheorem{corollary}[lemma]{Corollary}
\newtheorem*{claim}{Claim}
\newtheorem*{stat}{\name}
\newcommand{\name}{testing}

\theoremstyle{definition}
\newtheorem{definition}[lemma]{Definition}

\newtheorem{example}[lemma]{Example}
\newtheorem{problem}{Problem}
\newtheorem{remark}[lemma]{Remark}

\newcommand{\qedc}{{\qed}~{\rm Claim~{\theclaim}.}}

\numberwithin{equation}{section}

\newcommand{\ootimes}{\mathbin{\bar{\otimes}}}

\newcommand{\ootimest}{\mathbin{\vec{\otimes}}}

\newcommand{\jz}{$\set{\jj,0}$}

\newcommand{\M}{\E{M}}
\newcommand{\SD}{$(\mathrm{SD}_{\jj})$}

\DeclareMathOperator{\Int}{Int}
\DeclareMathOperator{\Cov}{Cov}

\DeclareMathOperator{\Max}{Max}
\DeclareMathOperator{\J}{J}
\DeclareMathOperator{\Hom}{Hom}
\DeclareMathOperator{\Id}{Id}
\DeclareMathOperator{\Con}{Con}
\DeclareMathOperator{\Conc}{Con_c}

\DeclareMathOperator{\Pow}{P} 

\newcommand{\CC}{\E{C}}
\newcommand{\T}{$(\mathrm{T})$}
\newcommand{\R}{\trianglelefteq}
\newcommand{\TJ}{$(\mathrm{T}_\jj)$}

\newcommand{\FL}{\mathrm{F}}

\newcommand{\fa}{{\F{a}}}

\newcommand{\dd}{\mathrm{d}}

\begin{document}

\title{Tensor products and transferability of semilattices}

\author{G.~Gr\"atzer}
\thanks{The research of the first author was partially
        supported by the NSERC of Canada.}
\address{Department of Mathematics\\
         University of Manitoba\\
         Winnipeg, Manitoba\\
         Canada R3T 2N2}
\email{gratzer@cc.umanitoba.ca}
\urladdr{http://www.maths.umanitoba.ca/homepages/gratzer.html/}

\author{F.~Wehrung}
\address{C.N.R.S., E.S.A. 6081\\
         D\'epartement de Math\'ematiques\\
         Universit\'e de Caen\\
         14032 Caen Cedex\\
         France}
\email{wehrung@math.unicaen.fr}
\urladdr{http://www.math.unicaen.fr/\~{}wehrung}

\date{March 25, 1999}

\begin{abstract}
In general, the tensor product, $A \otimes B$, of the lattices $A$ and
$B$ with zero is not a lattice (it is only a join-semilattice with
zero). If $A\otimes B$ is a \emph{capped} tensor product, then
$A\otimes B$ is a lattice (the converse is not known). In this paper, we
investigate lattices $A$ with zero enjoying the property that $A\otimes
B$ is a capped tensor product, for \emph{every} lattice $B$ with zero; we
shall call such lattices \emph{amenable}.

The first author introduced in 1966 the concept of a \emph{sharply
transferable lattice}.  In 1972, H.~Gaskill
\cite{Gask73} defined, similarly, sharply transferable semilattices,
and characterized them by a very effective condition \T.

We prove that \emph{a finite lattice $A$ is \emph{amenable} if{f} it is
\emph{sharply transferable} as a join-semilattice}.

For a general lattice $A$ with zero, we obtain the result: \emph{$A$ is
amenable if{f} $A$ is locally finite and every finite sublattice of $A$
is transferable as a join-semi\-lattice}.

This yields, for example, that a finite lattice $A$ is amenable
if{f} $A\otimes\FL(3)$ is a lattice if{f} $A$ satisfies \T, with
respect to $\jj$.  In particular,
$M_3\otimes\FL(3)$ is not a lattice. This solves a problem raised by
R.~W. Quackenbush in 1985 whether the tensor product of lattices with
zero is always a lattice.
\end{abstract}

\maketitle
\section{Introduction}

The tensor product, $A \otimes B$, of the \jz-semilattices $A$ and $B$
is defined in a very classical fashion, as a free (universal) object
with respect to a natural notion of \emph{bimorphism}, see G.~Fraser
\cite{Fras78}, G. Gr\"atzer, H. Lakser, and R.W. Quackenbush
\cite{GLQu81}, G. Gr\"atzer and F. Wehrung \cite{GrWe}.
Unfortunately, the tensor product of two lattices with zero is, in
general, not a lattice.

Let $A$ and $B$ be lattices with zero. If the tensor product, $A
\otimes B$, satisfies the very natural condition of being
\emph{capped}, introduced in
\cite{GrWe}, then $A \otimes B$ is always a lattice.  Capped tensor
products have many interesting properties.  The most important one is
the main result of \cite{GrWe}: If $A \otimes B$ is a capped tensor
product, then $\Conc(A \otimes B)$ and $\Conc A \otimes \Conc B$ are
isomorphic. For finite lattices this was proved in G. Gr\"atzer, H.
Lakser, and R.W. Quackenbush \cite{GLQu81}.

In this paper, we study the connections between $A\otimes B$ being a
lattice and $A\otimes B$ being capped, for lattices $A$ and $B$ with
zero. We do not have the answer to the most obvious question whether
these two conditions are equivalent (see Problem~\ref{Q:EquiSm} in
Section~\ref{S:problems}). We prove, however, that these two conditions
are equivalent provided that $A$ is finite (or locally finite), see
Theorem~\ref{T:SmFin}.

Since we cannot handle the general problem, when is $A \otimes B$ a
lattice, we universally quantify one of the variables: Let us call the
lattice $A$ with zero \emph{amenable}, if $A \otimes B$ is a capped
tensor product, for any lattice $B$ with zero. Then from
Theorem~\ref{T:SmFin} we obtain the easy corollary that for a
(locally) finite lattice $A$ with zero,
$A$ is amenable if{f} $A \otimes B$ is a lattice, for any lattice
$B$ with zero. This leads to the central problem of this paper:
characterize amenable lattices.

The technique to handle this problem comes from an unexpected source.

The first author introduced in 1966 the concept of \emph{sharp
transferability} for lattices (see \cite{gG70Ca}).  A finite lattice $A$
is sharply transferable, if whenever $A$ has an embedding $\gf$ into
$\Id L$, the ideal lattice of a lattice $L$, then $A$ has an embedding
$\gy$ into $L$ satisfying $\gy(x) \in \gf(y)$ if{f} $x \leq y$.  In 1972,
H.~Gaskill \cite{Gask73} characterized sharply transferable semilattices
by a very effective condition \T. For a lattice~$A$, as usual, we shall
denote by \TJ\ the condition \T\ for $\vv<A; \jj>$. Finite lattices
satisfying \TJ\ are well-understood. In particular, they are exactly the
finite, \emph{lower bounded} homomorphic images of finitely generated
free lattices, see R.~Freese, J.~Je\v zek, and J.B. Nation
\cite{FJNa95}, from which we shall borrow a number of results.

We characterize finite amenable lattices as follows:

\medskip

\emph{For a finite lattice $A$, the following conditions are equivalent:
 \begin{enumerate}
 \item $A$ is amenable.
 \item $A \otimes B$ is a lattice, for every lattice $B$ with zero.
 \item $A$ as a join-semilattice is sharply transferable.
 \item $A \otimes \FL(3)$ is a lattice.
 \item $A$ satisfies \TJ.
 \end{enumerate}}

\medskip

Our main result (see Corollary~\ref{C:SmFin} and Theorem~\ref{T:JoinAm})
characterizes amenability for general lattices:

\medskip

\emph{For a lattice $A$ with zero, the following conditions are
equivalent:
 \begin{enumerate}
\item $A$ is amenable.
\item $A$ is locally finite and
$A\otimes B$ is a lattice, for every lattice $B$ with zero.
\item $A$ is locally finite and $A \otimes \FL(3)$ is a lattice.
\item $A$ is locally finite and every finite sublattice of
$A$ satisfies \TJ.
\end{enumerate}}

\medskip

The finite case trivially follows from the general result.

Since $M_3$ fails \TJ, it follows that $M_3\otimes\FL(3)$ is not a
lattice. This answers a problem raised by R.W. Quackenbush in 1985
whether the tensor product of lattices with zero is always a lattice.
By a different approach, this was also answered in our paper
\cite{GrWe2}, where we produce two different examples of lattices $A$
and $B$ with zero such that
$A \otimes B$ is not a lattice: in the first example, $A$ and $B$ are
planar; in the second, $A$ and $B$ are modular.

Sections~\ref{S:antitone} and \ref{S:Transf} are introductory. In
Section~\ref{S:antitone}, we review the basic facts about tensor
products with special emphasis on the representation by some
hereditary subsets of
$A \times B$ and the representation by antitone maps, as in J.
Anderson and N. Kimura \cite{AnKi78}. Section~\ref{S:Transf}
introduces transferability and the technical tools we inherit from
transferability: minimal pairs, the condition \T, and the adjustment
sequence of a map. Section~\ref{S:Transf} introduces lower bounded
homomorphisms as well, together with some required technical tools,
such as the lower limit table.

In Section~\ref{S:capped}, we recall the definition of a capped tensor
product, and we establish a number of technical results; applying a
theorem of \cite{GrWe}, we prove, for instance, that if $A\otimes B$
is a capped tensor product, then so is $(A/\ga)\otimes(B/\gb)$, where
$\ga$ is a lattice congruence of $A$ and $\gb$ is a lattice congruence
of $B$.

Amenable lattices are introduced in Section~\ref{S:univcapped}, where we
characterize amenable locally finite lattices. In Section~\ref{S:Step},
we prepare the ground for studying capped tensor products \emph{via}
adjustment sequences of maps. Not all maps can be considered, but only
those that we call \emph{step functions}, which are ``measurable'' finite
joins of characteristic functions. In Section~\ref{S:ABSm}, we utilize the
concepts introduced in Sections~\ref{S:Transf} and~\ref{S:Step} to
characterize capped tensor products $A\otimes B$, with $A$ and $B$
arbitrary lattices with zero.

In Section~\ref{S:TJimplies}, we characterize amenable lattices, as
locally finite lattices with zero in which every finite sublattice
satisfies \TJ. Section~\ref{S:problems} concludes the paper
discussing some related results and stating some open problems.

\section{Tensor products}\label{S:antitone}
\subsection{The basic concepts}\label{S:basic} We shall adopt the
notation and terminology of our paper \cite{GrWe}. In particular, for
a \jz-semilattice $A$, we use the notation $A^-=A -
\set{0}$.  Note that $A^-$ is a subsemilattice of $A$.

Let $A$ and $B$ be \jz-semilattices. We denote by $A \otimes B$ the
\emph{tensor product} of $A$ and $B$, defined as the free
\jz-semilattice generated by the set $A^- \times B^-$ subject to the
relations $\vv<a,b_0>
\jj \vv<a,b_1> = \vv<a, b_0 \jj b_1>$, for $a  \in A^-$, $b_0$, $b_1
\in B^-$; and symmetrically, $\vv<a_0, b> \jj \vv<a_1, b> = \vv<a_0
\jj a_1, b>$, for $a_0$, $a_1  \in  A^-$,  $b  \in B^-$.

It follows directly from the definition that the sum distributes over
tensor product; since the sum for \jz-semilattices is the direct
product, we get the formula:
 \begin{equation}\label{E:distr}
 (A \times B) \otimes C \iso  (A \otimes C) \times (B \otimes  C).
 \end{equation}

The following two statements on tensor products are taken from
Corollary 3.7(iv) and Corollary 3.9 of \cite{GrWe}.  For a lattice
$L$, we denote by
$\Con L$ the congruence lattice of $L$.

\begin{lemma}\label{L:3.7}\hfill
 \begin{enumerate}
\item Let $A$ and $B$ be lattices with zero, let $\ga \in \Con A$ and
$\gb \in \Con B$. If $A \otimes B$ is a lattice, then
$(A/\ga) \otimes (B/\gb)$ is also a lattice.
\item  Let $A$, $A'$, $B$, $B'$ be lattices with zero such that $A$
is a sublattice of $A'$ and $B$ is a sublattice of $B'$. If
$A' \otimes B'$ is a lattice, then $A \otimes B$ is a lattice.
\end{enumerate}
 \end{lemma}

\subsection{The set representation}\label{S:set} In \cite{GrWe}, we
used the following representation of the tensor product.

First, we introduce the notation:
 \[
   \bot_{A,B} = (A \times \set{0})  \uu  (\set{0} \times B).
 \]

Second, we introduce a partial binary operation on $A \times B$: let
$\vv<a_0,b_0>$, $\vv<a_1,b_1> \in A \times B$; the \emph{lateral join}
of $\vv<a_0,b_0>$ and $\vv<a_1,b_1>$ is defined if $a_0 = a_1$ or
$b_0 = b_1$, in which case, it is the join, $\vv<a_0\jj a_1,b_0\jj b_1>$.

Third, we define bi-ideals: a nonempty subset $I$ of $A \times B$ is a
\emph{bi-ideal} of $A \times  B$, if it satisfies the following
conditions:
 \begin{enumerate}
 \item $I$ is hereditary;
 \item $I$ contains $\bot_{A,B}$;
 \item $I$ is closed under lateral joins.
 \end{enumerate}

The \emph{extended tensor product} of $A$ and $B$, denoted by $A
\ootimes B$, is the lattice of all bi-ideals of $A \times B$.

It is easy to see that $A \ootimes  B$ is an algebraic lattice. For
$a \in A$ and $b \in B$, we define $a \otimes b \in A \ootimes B$ by
 \[
   a \otimes b =\bot_{A,B}  \uu  \setm{\vv<x, y>  \in A \times B}
  {\vv<x, y>  \leq \vv<a, b>}
 \]
 and call $a \otimes b$ a \emph{pure tensor}.  A pure tensor is a
principal (that is, one-generated) bi-ideal.  Using this notation, the
isomorphism \eqref{E:distr} can be established by the map
 \begin{equation}\label{E:map}
   \vv<a, b> \otimes c \mapsto \vv<a \otimes c, b \otimes c>.
 \end{equation}

Now we can state the representation:

\medskip

\emph{The tensor product $A \otimes B$ can be represented as the
\jz-sub\-semi\-lat\-tice of compact elements of $A \ootimes B$.}

\medskip

For every positive integer $m$, denote by $\FL(m)$ the free lattice on
$m$ generators $x_0$, \dots, $x_{m - 1}$. One can evaluate any element
$p$ of $\FL(m)$ at any $m$-tuple of elements of any lattice, thus
giving the notation $p(a_0,\ldots,a_{m-1})$. Let $p^\dd$ be the dual
of $p$.

The following purely arithmetical formulas are due to G.~A. Fraser
\cite{Fras78}.  They are easiest to prove using the above
representation of tensor products.

\begin{lemma}\label{L:formulas}
Let $A$ and $B$ be \jz-semilattices.
Let $a_0$, $a_1 \in A$, and $b_0$,
$b_1  \in B$ such that $a_0 \mm  a_1$ and $b_0 \mm  b_1$ both exist.
\begin{enumerate}
\item The intersection (meet) of two pure tensors is a pure tensor, in
fact,
 \begin{equation*}
   (a_0 \otimes b_0) \ii (a_1 \otimes b_1)
     = (a_0 \mm  a_1) \otimes (b_0 \mm  b_1).
 \end{equation*}
 \item The join of two pure tensors is the union of four pure tensors,
in fact,
 \begin{align*}
      (a_0 \otimes b_0) \jj (a_1 \otimes b_1) & = \\
    (a_0 \otimes b_0)  \uu  (a_1 \otimes b_1)
     & \uu  ((a_0 \jj a_1) \otimes (b_0 \mm  b_1))
  \uu  ((a_0 \mm  a_1) \otimes (b_0 \jj b_1)).
 \end{align*}
 \item Let $A$ and $B$ be lattices with zero, let $n$ be a positive
integer, let $a_0$, \dots, $a_{n - 1}  \in A$, and let
$b_0$, \dots, $b_{n - 1}  \in B$. Then
 \[
	\JJm{a_i \otimes b_i}{i < n} = \UUm{p(a_0, \ldots, a_{n - 1}) \otimes
 p^\dd(b_0, \ldots, b_{n - 1})}{p  \in \FL(n)}.
 \]
 \item Let $A$ and $B$ be lattices with zero.  Then
 \begin{gather*}
  \JJm{a_i \otimes b_i}{i < n} \mm \JJm{c_j \otimes d_j}{j < m}\\
     = \UU (p(a_0, \dots, a_{n - 1}) \mm q(c_0, \dots, c_{m - 1}))
      \otimes (p^{\mrm{d}}(b_0, \dots, b_{n - 1}) \mm
      q^{\mrm{d}}(d_0, \dots, d_{m - 1})),
 \end{gather*}
where $p^{\mrm{d}}$ and $q^{\mrm{d}}$ are the duals of
$p$ and $q$, respectively, and where the union is for all $p  \in
\FL(n)$ and $q \in \FL(m)$.
 \end{enumerate}
 \end{lemma}

Note that in (iii) (and similarly, in (iv)) the right side is an infinite
union, for $n > 2$.

For any subset $X$ of a lattice $L$,  let $X^\mm$ (resp., $X^\jj$)
denote the meet-subsemilattice (resp., join-subsemilattice) of $L$
generated by $X$. We also write $X^{\mm\jj}=(X^\mm)^\jj$,
$X^{\jj\mm}=(X^\jj)^\mm$, and so on.
We define inductively an increasing sequence of finite
\emph{join-subsemilattices} $S_n(m)$ of $\FL(m)$ as follows:
 \begin{enumerate}
 \item $S_0(m)=\set{x_0, \ldots, x_{m-1}}^{\mm\jj}$.
 \item $S_{n + 1}(m)=S_n(m)^{\mm\jj}$.
 \end{enumerate}
In particular, all the $S_n(m)$ are finite and their
union equals $\FL(m)$.

The following lemma readily follows from Lemma~\ref{L:formulas}(iii).

\begin{lemma}\label{L:TensFinUn}
 Let $A$ and $B$ be lattices with zero, let $m$ be a positive integer,
let $a_0$, \dots, $a_{m - 1} \in A$, and let $b_0$, \dots,
$b_{m - 1} \in B$. Then the following conditions are equivalent:
 \begin{enumerate}
 \item The finite \emph{join} $\JJm{a_i\otimes b_i}{i<m}$ is a finite
\emph{union} of pure tensors.
 \item There exists a positive integer $n$ such that the following
finite \emph{union}:
 \[
   \UUm{p(a_0,\ldots,a_{m-1})\otimes p^\dd(b_0,\ldots,b_{m-1})}
   {p\in S_n(m)}
 \]
 belongs to $A\otimes B$.
 \end{enumerate}
 \end{lemma}

\subsection{Representation by homomorphisms}\label{S:homomorphisms}
Let $A$ and $B$ be \jz-semilattices. Note that $\Id B$, the set of all
ideals of $\vv<B; \jj>$, is a semilattice under intersection. So we can
consider the set of all semilattice homomorphisms from the semilattice
$\vv<A^-; \jj>$ into the semilattice $\vv<\Id B; \ii>$,
 \[
    A \ootimest B = \Hom(\vv<A^-; \jj>, \vv<\Id B; \ii>),
 \] ordered componentwise, that is, $f\leq g$ if{f} $f(a) \leq g(a)$
(that is,
$f(a) \ci g(a)$), for all $a \in A^-$. The arrow indicates which way
the homomorphisms go. Note that the elements of $A\ootimest B$ are
\emph{antitone} functions from $A^-$ to $\Id B$; indeed, if
$a \leq b$, then $a \jj b = b$ and so $f(a) \ii f(b) = f(b)$, that is,
$f(a) \ce f(b)$.

With any element $\gf$ of $A\ootimest B$, we associate the subset
$\ge(\gf)$ of $A \times B$:
 \[
  	\ge(\gf) = \setm{\vv<x, y> \in A \times B}{y \in \gf(x)} \uu
\bot_{A, B}.
 \]

\begin{proposition}\label{P:TensAnti} The map $\ge$ is an isomorphism
between $A\ootimest B$ and $A\ootimes B$. The~inverse map, $\ge^{-1}$,
sends $H \in A\ootimes B$ to $\ge^{-1}(H)
\colon A^- \to \Id B$, defined by
 \[
   \ge^{-1}(H)(a) = \setm{x  \in B}{\vv<a, x> \in H}.
 \]
 \end{proposition}

\begin{proof}
 We leave the easy computation to the reader.
\end{proof}

For $a  \in A$ and $b  \in B$, we can describe $\gx = \ge^{-1}(a
\otimes b)
\colon A^-\to\Id B$ as follows
 \[
 \gx(x)=
   \begin{cases}
    (b], &\text{if $x\leq a$};\\
   \set{0}, &\text{otherwise.}
   \end{cases}
 \]

If $A$ is \emph{finite}, then a homomorphism from $\vv<A^-;\jj>$ to
$\vv<\Id B;\ii>$ is determined by its restriction to $\J(A)$, the set
of all join-irreducible elements of $A$.

This representation of the tensor product is also utilized in
\cite{GrWe2} and \cite{GrWe3}.

\section{Sharply transferable lattices}\label{S:Transf}

In this section, we introduce sharply transferable lattices and
semilattices, minimal pairs, the condition \T, the $D$-relation, and the
adjustment sequence of a map. For more information on these topics, the
reader is referred to \cite{FJNa95} and \cite{GGPl75}.

\subsection{Sharp transferability}\label{S:Transferability}
We start with the definition of sharp transferability (see \cite{gG70Ca}).
We get a definition that is easier to utilize by using choice
functions ($\Pow(Y)$ is the power set of $Y$):

\begin{definition}
Let $X$ and $Y$ be sets and let $\gf\colon X\to \Pow(Y)$ be a map.
Then a \emph{choice function} for $\gf$ is a map $\gx\colon X\to Y$
such that $\gx(x)\in\gf(x)$, for all $x\in X$.
\end{definition}

\begin{definition}\label{D:sharplytransferable}
 A lattice $S$ is \emph{sharply transferable}, if for every embedding
$\gf$ of $S$ into $\Id T$, the ideal lattice  of a lattice $T$, there
exists an embedding $\gx$ of $S$ into~$T$ such that $\gx$ is a choice
function for
$\gf$ satisfying $\gx(x) \in \gf(y)$ if{f} $x \leq y$.
\end{definition}

Equivalently, $\gx(x) \in \gf(x)$ and $\gx(x) \nin \gf(y)$, for any $y
< x$.

The motivation for these definitions comes from the fact that the
well-known result: \emph{a~lattice $L$ is modular if{f} $\Id L$ is
modular}, can be recast: \emph{$N_5$ is a (sharply) transferable lattice}.

The reader should find it easy to verify that $N_5$ is a sharply
transferable lattice.  It is somewhat more difficult to see the
negative result: $M_3$ is not a sharply transferable lattice.

\subsection{Minimal pairs}\label{S:MinimalPairs}

Let $P$ be a poset and let $X$ and $Y$ be subsets of $P$. Then the
relation $X$ \emph{is dominated by} $Y$, in notation, $X \ll Y$, is
defined as follows:

\medskip

$X \ll Y$ if{f} for all $x \in X$ there exists $y \in Y$ such that $x
\leq y$.

\medskip

\noindent The relation $\ll$ defines a \emph{quasi-ordering} on the
power set of $P$, and it is a \emph{partial ordering} on the set of
antichains of
$P$.

The following definition is from H. Gaskill \cite{Gask73}:

\begin{definition}\label{D:MinPair} Let $A$ be a finite
join-semilattice. A \emph{minimal pair} of $A$ is a pair $\vv<p, I>$
such that $p \in \J(A)$, $I \ci \J(A)$, and the following three
conditions hold:
\begin{enumerate}
\item $p \nin I$.
\item $p \leq \JJ{I}$.
\item For all $J \ci \J(A)$, if $J \ll I$ and $p \leq \JJ{J}$, then
$I \ci J$.
\end{enumerate}
\end{definition}

If $\vv<p, I>$ is a minimal pair, then $I$ is an antichain of $\J(A)$.
If $I$ is an antichain of $\J(A)$ and $\vv<p, I>$ satisfies conditions
\ref{D:MinPair}(i), \ref{D:MinPair}(ii), then condition
\ref{D:MinPair}(iii) is equivalent to the following condition:
\begin{itemize}
\item[\tup{(iii')}] For all $J\ci\J(A)$, if $J\ll I$ and
$p\leq\JJ{J}$, then $I\ll J$.
\end{itemize}

``Minimal pairs'' are minimal in two ways.  Firstly, if $\vv<p, I>$ is
a minimal pair, then $p \leq \JJ{I}$ and an element of $I$ cannot be
replaced by a set of smaller join-irreducible elements while retaining
that their join is over $p$.  Indeed, if $i \in I$ can be replaced by
$X \ci \J(A)$ with $X \ii I = \es$, that is, $p \leq \JJ J$, where $J
= (I - \set{i}) \uu X$, then $J \ll I$, hence by (iii), $I\ci J$, a
contradiction.  Secondly, the collection of all minimal pairs is the
minimal information necessary to describe the join structure of the
finite lattice.

If $A$ is a finite join-semilattice and $p \in\J(A)$, then we shall define
 \[
 \M(p)=\setm{I \ci \J(A)}{\vv<p, I>\text{ is a minimal pair}}.
 \]

\subsection{The condition \T}\label{S:condT} For the following
definition, see H. Gaskill \cite{Gask73}.

\begin{definition} A finite join-semilattice $A$ satisfies \T, if
$\J(A)$ has a linear order $\R$ such that for every minimal pair
$\vv<p,J>$ of $A$, the relation $p \R j$ holds, for any $j \in J$.
\end{definition}

A finite lattice $A$ satisfies the condition \TJ, if the semilattice
$\vv<A;\jj>$ satisfies~\T.

One can define sharply transferable semilattices by changing
``lattice'' to ``semilattice'' in
Definition~\ref{D:sharplytransferable}. H. Gaskill
\cite{Gask73} proved that \T\ characterizes finite sharply transferable
semilattices.

An interpretation of \TJ\ for finite lattices can be found in
\cite{FJNa95}. Let $A$ be a finite semilattice. For $p\in\J(A)$, let
$p_*$ denote the unique element covered by $p$. The \emph{dependency
relation} is the binary relation $D$ defined on $\J(A)$ as follows

\medskip

$p\mathrel{D}q$\q if{f}\q $p \ne q$ and there exists $x \in A$ such that
$p\leq q\jj x$ and $p\nleq q_*\jj x$.

\medskip

\begin{lemma}\label{L:DMin}
Let $A$ be a finite lattice and $p$, $q\in\J(A)$. Then
$p\mathrel{D}q$ if{f} there exists
$I\in\M(p)$ such that $q \in I$.
\end{lemma}

\begin{corollary}\label{C:DMin}
A finite lattice $L$ satisfies \TJ\ if{f} $\J(L)$ has no $D$-cycle.
\end{corollary}

\subsection{Lower bounded homomorphisms; the lower limit table}
\label{S:LowBded}
We recall in this section some classical concepts about lower bounded
homomorphisms, due to A. Day, H.S. Gaskill, B. J\'onsson,
A. Kostinsky, and  R.N. McKenzie, which are presented in Chapter~2 of
R. Freese, J. Je\v{z}ek, and J.B. Nation \cite{FJNa95}.

Let $K$ and $L$
be lattices. If $h\colon K\to L$ is a lattice homomorphism, we say that
$h$ is \emph{lower bounded}, if the set $h^{-1}[a)$ has a least element,
for all $a\in L$.

Let us assume that $K$ is finitely generated and let $X$ be a finite
generating set of $K$. We put $H_0=X^{\mm}\uu\set{1_K}$, and
$H_n=X^{\mm(\jj\mm)^n}$, for all $n>0$.
We note that $\vv<H_n \mid n\in\go>$ is an increasing sequence
of \emph{finite} meet-subsemilattices of $K$, and that their union is
$K$. For $n\in\go$ and $a\in L$, we denote by $\gb_n(a)$ the least
element $x$ of $H_n$ such that $a\leq x$, if it exists. Note that, for
$n\in\go$, $\gb_n(a)$ is defined if{f} $a\leq h(1_K)$. If this holds,
then $\gb_n(a)$ can be computed by the following formulas, see
Theorem~2.3 in \cite{FJNa95} ($h[X]$ denotes the image of $X$ under $h$):
 \begin{align}
 \gb_0(a)&=\MMm{x\in X}{a\leq h(x)},
 \label{Eq:beta0}\\
 \gb_{n+1}(a)&=\gb_n(a)\mm\MM\left(\,\JJ\gb_n[S]\mid
 S\in\CC^*(a),\,\JJ S\leq h(1_K)\,\right),\label{Eq:betan}
 \end{align}
where an empty meet in $K$ equals $1_K$ and
 \[
 \CC^*(a)=\left\{\,S\ci K
 \mid S\text{ is finite nonempty},\,x\leq\JJ S,
 \text{ and }x\nin(S]\,\right\},
 \]
for all $a\in K$.

The sequence $\vv<\gb_n \mid n\in\go>$ is called the \emph{lower limit
table} of $h$. It is uniquely determined by $h$ and $X$.
The homomorphism $h$ is lower bounded if{f}, for all $a\leq h(1_K)$,
there exists $n\in\go$ such that $\gb_n(a)=\gb_{n+1}(a)$. In that
case, $\gb_n(a)$ is the least element of $h^{-1}[a)$.

A finitely generated lattice $L$ is \emph{lower bounded}, if there
exists a \emph{surjective}, lower bounded lattice homomorphism from a
finitely generated free lattice onto $L$. We recall the following
characterization of lower bounded lattices, see Theorem~2.13
in~\cite{FJNa95}:

\begin{proposition}\label{P:CharLB}
For a finitely generated lattice $L$, the following are equivalent:
\begin{enumerate}
\item $L$ is lower bounded.

\item For every finitely generated lattice $K$, every lattice
homomorphism from $K$ to $L$ is lower bounded.
\end{enumerate}
\end{proposition}

Finite lower bounded lattices have a number of characterizations.
The following result states one of them. For a proof, we
refer to Corollary~2.39 in \cite{FJNa95}.

\begin{proposition}\label{P:LBfin}
A finite lattice is lower bounded if{f} it satisfies \TJ.
\end{proposition}

We recall that another characterization of \TJ\ is given by
Corollary~\ref{C:DMin}.

We shall also need the following result, see the proof of Theorem~2.6
in \cite{FJNa95}:

\begin{proposition}\label{P:LoBdn}
Let $K$ be a finitely generated lattice, let $L$ be a lattice, and let
$h\colon K\twoheadrightarrow L$ be a \emph{surjective} lattice
homomorphism. Let
$X$ be a finite generating subset of $K$, and let
$\vv<\gb_n \mid n\in\go>$ be the associated lower limit table of $h$.
Let $n\in\go$, and let us assume that $\gb_n=\gb_{n+1}$. Then
$L=h[X]^{\mm(\jj\mm)^n}$, so that $L$ is finite. Furthermore, $h$ is
lower bounded.
\end{proposition}

The proof of of Proposition~\ref{P:LoBdn} is easy to outline.
We note first that $\gb_n(a)\in X^{\mm(\jj\mm)^n}$, for all $a\in L$.
Since $\gb_n=\gb_{n+1}$, it follows that $\gb_n(a)$ is the least element
of $h^{-1}[a)$, for all $a\in A^-$. In particular, $h$ is lower bounded.
Since $h$ is surjective, the equality
$a=h(\gb_n(a))$ holds. It follows that $a\in h[X]^{\mm(\jj\mm)^n}$.

\section{Capped tensor products of lattices}\label{S:capped}

In \cite{GrWe} (see Definition 7.1), we introduced \emph{capped tensor
products}. We recall the definition here:

\begin{definition}\label{D:cappedpair}
 Let $A$ and $B$ be lattices with zero. We say that $A\otimes B$ is a
\emph{capped tensor product}, if every element of $A\otimes B$ is a
\emph{finite union} of pure tensors.
 \end{definition}

It is easy to see by Lemma~\ref{L:formulas} that if $A\otimes B$ is a
capped tensor product, then $A\otimes B$ is a lattice (see Lemma 7.2 of
\cite{GrWe}).

In this section, we shall establish some basic properties of capped
tensor products. Most of the results in this section are technical
lemmas with the exception of Theorem~\ref{T:product}, stating that
$(A\times B)\otimes C$ is capped if{f} both $A\otimes C$  and
$B\otimes C$ are capped.

\begin{lemma}\label{L:aibiTens}
 Let $A$ and $B$ be lattices with zero. Let $n$ be a positive integer,
let
$a_0$, \dots, $a_{n - 1}  \in A$, and let $b_0$, \dots, $b_{n - 1}
\in B$. Then
 \[
   \UUm{a_i \otimes b_i}{i < n}  \in A \otimes B
 \]
if{f}, for all $i$, $j < n$, $i \neq j$, the following two
conditions hold:
 \begin{enumerate}
 \item $a_i \mm a_j = 0$ or there exists $k < n$ such that  $a_i \mm
a_j \leq a_k$ and $b_i \jj b_j \leq b_k$;
 \item $b_i \mm b_j = 0$ or there exists $k < n$ such that  $b_i \mm
b_j \leq b_k$ and $a_i \jj a_j \leq a_k$.
 \end{enumerate}
 \end{lemma}

\begin{proof}
 It is evident that the hereditary subset $H = \UUm{a_i \otimes
b_i}{i<n} \ci A \times B$ belongs to $A \otimes B$ if{f} it is closed
under lateral joins. The two conditions guarantee this. Now the
conclusion easily follows.
 \end{proof}

\begin{corollary}\label{C:aibiTens}
 Let $A$, $A'$, $B$, and $B'$ be lattices with zero, let
$f \colon A\to A'$ and $g \colon B\to B'$ be lattice homomorphisms.
Let $n$ be a positive integer, let
$a_0$, $\ldots$, $a_{n-1} \in A$, and let $b_0$, $\ldots$,
$b_{n-1}\in B$.
 \begin{itemize}
 \item[\tup{(a)}] If $f$ and $g$ are $0$-preserving, then $\UUm{a_i
\otimes b_i}{i < n} \in A\otimes B$ implies that $\UUm{f(a_i) \otimes
g(b_i)}{i < n} \in A' \otimes B'$.
 \item[\tup{(b)}] If $f$ and $g$ are lattice embeddings, then
$\UUm{f(a_i)
\otimes g(b_i)}{i < n} \in A' \otimes B'$ implies that $\UUm{a_i
\otimes b_i}{i < n}  \in A \otimes B$.
 \end{itemize}
 \end{corollary}

\begin{corollary}\label{C:aibiFinUn}
Let $A$, $A'$, $B$, and $B'$ be
lattices with zero, let $f \colon A \to A'$ and $g \colon B \to B'$ be
lattice homomorphisms. Let $n$ be a positive integer, let $a_0$,
\dots, $a_{n - 1} \in A$ and let $b_0$, \dots, $b_{n - 1} \in B$.
\begin{enumerate}
\item Let $f$ and $g$ be $0$-preserving. If $\JJm{a_i\otimes
b_i}{i<n}$ is a finite union of pure tensors in $A\otimes B$, then
$\JJm{f(a_i)\otimes g(b_i)}{i<n}$ is a finite union of pure tensors
in $A'\otimes B'$.
\item Let $f$ and $g$ be lattice embeddings. If $\JJm{f(a_i)\otimes
g(b_i)}{i<n}$ is a finite union of pure tensors in $A'\otimes B'$,
then $\JJm{a_i\otimes b_i}{i<n}$ is a finite union of pure tensors in
$A\otimes B$.
\end{enumerate}
\end{corollary}

Now we are ready to state our first preservation results about capped
tensor products:

\begin{proposition}\label{P:SmSubQuot} Let $A\otimes B$ be a capped
tensor product of the lattices $A$ and $B$ with zero. Then the
following holds:
\begin{enumerate}
\item Let $\ga\in\Con A$ and $\gb\in\Con B$. Then
$(A/\ga)\otimes(B/\gb)$ is a capped tensor product of $A/\ga$ and
$B/\gb$.
\item Let $A'\ci A$ and $B'\ci B$ be lattices with zero. Then
$A'\otimes B'$ is a capped tensor product of $A'$ and $B'$.
\end{enumerate}
\end{proposition}

Note that we do \emph{not} assume in (ii) that $0_A = 0_{A'}$ or that
$0_B = 0_{B'}$.

\begin{proof} (i) follows immediately from
Corollary~\ref{C:aibiFinUn}(i), while (ii) follows from
Corollary~\ref{C:aibiFinUn}(ii).
\end{proof}

\begin{proposition}\label{P:SmLim}
Let $A$ and $B$ be
lattices with zero, and let $(A_i \mid i \in I)$ be a directed
family of lattices with zero. Suppose that, with appropriate
transition maps that are $0$-lattice homomorphisms, $A = \varinjlim_i
A_i$, and that all $A_i\otimes B$ are capped tensor products. Then
$A\otimes B$ is a capped tensor product.
\end{proposition}

Now we present the result on direct products:

\begin{theorem}\label{T:product} Let $A$, $B$, and $C$ be lattices with
zero. Then $(A\times B)\otimes C$ is a capped tensor product if{f} both
$A\otimes C$ and $B\otimes C$ are capped tensor products.
\end{theorem}

\begin{proof} If $(A\times B)\otimes C$ is a capped tensor product,
then both $A\otimes C$ and $B\otimes C$ are capped tensor products by
Proposition~\ref{P:SmSubQuot}(i).

Conversely, let us assume that both $A\otimes C$ and $B\otimes C$ are
capped tensor products. Let $H \in (A \times B) \otimes C$. We have to
prove that $H$ is a finite union of pure tensors. By \eqref{E:distr},
we take instead $H \in (A \otimes C) \times (B \otimes C)$.

Let $\vv<u, v> \in (A \otimes C) \times (B \otimes C)$; then $u$ is a
finite union of pure tensors $a \otimes c \in A \otimes C$ and $v$ is
a finite union of pure tensors $b \otimes c \in B \otimes C$. So
$\vv<u, v>$ is a finite union of pure tensors of the form $\vv<a
\otimes c_1,b \otimes c_2> \in (A \otimes C) \times (B \otimes C)$.
Since
\[
\vv<a \otimes c_1, b \otimes c_2> = \vv<a \otimes c_1, 0> \jj \vv<0, b
\otimes c_2>,
\] by formula (ii) of Lemma~\ref{L:formulas}, $\vv<a \otimes c_1, b
\otimes c_2>$ is a union of (at most) four pure tensors; therefore,
$(A\times B)\otimes C$ is a capped tensor product.
\end{proof}

By using Proposition~\ref{P:SmLim}, we deduce immediately the
following:

\begin{corollary}\label{C:SmInfOplus} Let $(A_i \mid i \in I)$ be a
family of lattices with zero and let $A$ be the discrete direct
product of this family. Then, for any lattice $B$ with zero, if
$A_i\otimes B$ is a capped tensor product, for all $i \in I$, then
$A\otimes B$ is a capped tensor product.
\end{corollary}

\section{Amenable lattices}\label{S:univcapped}
Now we come to the central concept of this paper:

\begin{definition}\label{D:SmLatt}
A lattice $A$ with zero is \emph{amenable}, if $A\otimes L$ is a
capped tensor product, for every lattice $L$ with zero.
\end{definition}

In other words, $A$ is amenable if{f} $A \otimes L$ is always
capped; it follows that then $A \otimes L$ is always a
lattice.
\begin{theorem}\label{T:SmPres}
The class of amenable
lattices with zero is preserved under the following operations:
\begin{enumerate}
\item the formation of sublattices,
\item the formation of quotient lattices,
\item finite direct products,
\item direct limits.
\end{enumerate}
\end{theorem}

It follows from known results on lower bounded lattices that
Theorem~\ref{T:JoinAm} is stronger than Theorem~\ref{T:SmPres}. However,
Theorem~\ref{T:SmPres} is more elementary, and it can easily be
generalized to $\mbf{V}$-amenable lattices, for any variety $\mbf{V}$ of
lattices, see Definition~\ref{D:CJA}.

We shall see, in Corollary~\ref{C:SmNonVar}, that the class of amenable
lattices is \emph{not} closed under arbitrary direct products; in
particular, it is not a variety. However, the following result holds:

\begin{corollary}\label{C:SmFinVar}
Let $\mbf V$ be a variety
generated by a amenable finite lattice. Then every lattice with
zero in $\mbf V$ is amenable.
\end{corollary}

\begin{proof}
Let $A$ be a amenable finite lattice generating
$\mbf V$. Let $B \in \mbf V$ and let $U$ be a finitely generated, say,
$n$-generated, $\set{0}$-sublattice of $B$. Let $\FL_{\mbf{V}}(n)$ be the
free lattice on $n$ generators in $\mbf{V}$. Since $U$ belongs to
$\mbf{V}$ and $U$ is $n$-generated, there exists a surjective lattice
homomorphism $\pi \colon \FL_{\mbf{V}}(n) \twoheadrightarrow U$. But
$\FL_{\mbf{V}}(n)$ embeds into $A^{A^n}$ (this is a classical result
of universal algebra). Since $A$ is amenable and $A^n$ is finite,
it follows from Theorem~\ref{T:SmPres} that $\FL_{\mbf{V}}(n)$ is
amenable. Since $U$ is a quotient of $\FL_{\mbf{V}}(n)$, $U$ is amenable,
again by Theorem~\ref{T:SmPres}.

Finally, $B$ is the direct union of all its finitely generated
$\set{0}$-sublattices, so the conclusion follows by
Proposition~\ref{P:SmLim}.
\end{proof}

Again, Corollary~\ref{C:SmFinVar} can easily be
generalized to $\mbf{W}$-amenable lattices, for any variety
$\mbf{W}$ of lattices.

In Example~\ref{E:LocFin}, we show that
``amenable finite lattice''
cannot be replaced by ``amenable locally finite lattice'' in
Corollary~\ref{C:SmFinVar}.

\begin{corollary}\label{C:DistrSm}
Every distributive lattice with zero is amenable.
\end{corollary}

If $A$ is an amenable lattice with zero, then $A \otimes L$ is a
lattice, for every lattice $L$ with zero. We do not know whether the
converse is true in general (see Problem~\ref{Q:EquiSm} in
Section~\ref{S:problems}), however, we can settle this problem for
locally finite lattices:

\begin{theorem}\label{T:SmFin}
Let $A$ and $B$ be lattices with zero.
Let $A$ be \emph{locally finite}. Then $A \otimes B$ is a lattice if{f}
$A\otimes B$ is a capped tensor product.
\end{theorem}

\begin{proof} We prove the nontrivial direction; so we assume that $A
\otimes B$ is a lattice.

\emph{First case: $A$ is finite.} Let $H \in A \otimes B$; we have to
prove that $H$ is a finite union of pure tensors. Let $\gx \colon
\vv<A^-; \jj> \to \vv<\Id B; \ii>$ be the antitone map associated
with~$H$, as defined in Section~\ref{S:homomorphisms}. Then $H = \UUm{a
\otimes \gx(a)}{a \in A^-}$ (where $a \otimes \gx(a)$ is the set of all
$a \otimes x$, $x \in \gx(a)$), so it suffices to prove that all
$\gx(a)$, $a \in A^-$, are principal ideals of $B$.

For all $a \in A^-$,
\[
\gx(a)=\IIm{\gx(p)}{p \in \J(A) \text{ and }p\leq a}; \] thus it
suffices to prove that $\gx(p)$ is a principal ideal of $B$, for every
$p \in \J(A)$.

Let $u_B \in B$ such that $H\ci 1_A\otimes u_B$; this element $u_B$
exists because $H$ is compact.

In Proposition 4.3 of \cite{GrWe}, we noted that if a tensor product
is a lattice, then the set representation is closed under
intersection. Hence
\begin{align*} U &= H \ii (p_* \otimes u_B),\\ V &= H\ii(p \otimes u_B)
\end{align*} belong to $A\otimes B$. Since $U \ci V$, there exists a
decomposition of $V$ of the form
\begin{equation}\label{Eq:V/U} V = U \jj \JJm{a_i \otimes b_i}{i < n},
\end{equation} where $n$ is a nonnegative integer, and $0_A < a_i$ and
$0_B < b_i\leq u_B$, for all $i < n$. The~inequality $a_i \otimes b_i \leq
V$ holds in $A \otimes B$, for all $i < n$. Since $a_i$ and $b_i$ are
nonzero, $a_i\leq p$, therefore, $a_i\leq p_*$ or $a_i = p$. But if $a_i
< p$, then $a_i \otimes b_i \leq U$, by the definition of $U$, so that
$a_i \otimes b_i$ is absorbed by $U$ in the decomposition (\ref{Eq:V/U}).
In~other words, in (\ref{Eq:V/U}) we may assume that $a_i = p$, for all
$i<n$.

Set $b = \JJm{b_i}{i<n}$; then
 \[
 V = U\jj(p \otimes b),
 \]
and we define $V' = U \uu(p \otimes b)$. We claim that $V'$ is a
bi-ideal. It is obvious that $V'$ is a hereditary subset of $A
\otimes B$ containing $\bot_{A,B}$. We show that $V'$ is closed
under lateral joins.

Let $\vv<x_0, y>$ and $\vv<x_1, y>$ in $V'$. Without loss of
generality, we can assume that $x_0$, $x_1$, and $y$ are nonzero, and
$\vv<x_0, y> \in U$, $\vv<x_1, y> \in p \otimes b$. Then $x_0 \jj
x_1\leq p$ and $y\leq b$, so that $\vv<x_0 \jj x_1,y> \in V'$.

Let $\vv<x,y_0>$ and $\vv<x,y_1>$ in $V'$. Without loss of
generality, $x$, $y_0$, and $y_1$ are nonzero, and $\vv<x, y_0> \in
U$, $\vv<x,y_1> \in p \otimes b$. Since $U \leq p_* \otimes u_B$, it
follows that $x \leq p_*$, so
 \[
 \vv<x, y_0 \jj y_1> \in H \ii(p_* \otimes u_B) = U \leq V'.
 \]
This proves that $V'$ belongs to $A \otimes B$, hence $V = V'$.
Therefore,
 \[
 V = U\uu(p \otimes b).
 \]
It follows that, for all $y \in B$,
 \begin{align*}
 y \in\gx(p)&\text{ if{f} }p \otimes y \ci H,\\
 &\text{ if{f} }p \otimes y \ci V,\\
 &\text{ if{f} }p \otimes y \in U\uu(p \otimes b),\\
 &\text{ if{f} }y \leq b.
 \end{align*}
This proves that $\gx(p) = (b]$, a principal ideal, thus
completing the proof of the finite~case.

\emph{Second case: $A$ is locally finite.} Then $A$ can be written as
the direct union of all of its finite $\set{0}$-sublattices. By
Lemma~\ref{L:3.7}, the tensor product $U \otimes B$ is a lattice, for
each of those lattices~$U$. Thus by the first case, $U\otimes B$ is a
capped tensor product. Therefore, by Proposition~\ref{P:SmLim},
$A\otimes B$ is a capped tensor product.
\end{proof}

\begin{corollary}\label{C:SmFin} Let $A$ be a locally finite lattice
with zero. Then the following conditions are equivalent:
\begin{enumerate}
\item $A$ is amenable.
\item $A \otimes L$ is a lattice, for every lattice $L$ with zero.
\item $A \otimes\FL(3)$ is a lattice.
\end{enumerate}
\end{corollary}

\begin{proof} The equivalence of (i) and (ii) follows immediately from
Theorem~\ref{T:SmFin}. (ii)~implies (iii) is trivial.

Finally, assume (iii). For every set $X$, denote by $\FL^{(0)}(X)$ the
free $\set{0}$-lattice on the generating set $X$. By a result of P. M.
Whitman \cite{Whit41} (see Theorem VI.2.8 in~\cite{Grat98}),
$\FL^{(0)}(\go)$ embeds into $\FL(3)$ as a $\set{0}$-sublattice ($\go$
denotes the set of all natural numbers). Thus, by Lemma
\ref{L:3.7}(ii), $A \otimes\FL^{(0)}(\go)$ is a lattice. For every
infinite set~$X$, $\FL^{(0)}(X)$ is a direct limit of $0$-lattices
each one isomorphic to $\FL^{(0)}(\go)$; thus $A \otimes
\FL^{(0)}(X)$ is a lattice. But every lattice with zero is a quotient
of some $\FL^{(0)}(X)$, so we conclude the argument by
Lemma~\ref{L:3.7}(i).
\end{proof}

We have not yet seen a finite lattice that is not amenable. The
smallest example is $M_3$, as we shall see it in the next few
sections.

\section{Step functions and adjustment sequences}\label{S:Step}
The referee informed us that some of the arguments in this section date back
to A. Kostinsky's unpublished 1971 thesis, which is unavailable to us.

Let $A$ and $B$ be lattices with zero. We shall put
$(a]_{\bullet}=(a]\ii A^-$, for all $a\in A$. We shall denote by $\Int A$
the Boolean lattice of $\Pow(A^-)$ generated by all sets of the form
$(a]_{\bullet}$, for $a\in A^-$. Furthermore, we shall denote by
$\Int_*A$ the ideal of $\Int A$ consisting of all $X \in \Int A$
with $X\ci(a]_{\bullet}$, for some $a\in A$. For every subset $U$
of $A^-$, we denote by $\Max U$ the set of all \emph{maximal} elements of
$U$.

\subsection{Step functions}

We start with the following very simple lemma:

\begin{lemma}\label{L:MaxInt}
Let $U\in\Int_*A$. Then the two following properties hold:
 \begin{enumerate}
 \item Every element of $U$ is contained in an element of $\Max U$.
 \item $\Max U$ is finite.
 \end{enumerate}
 \end{lemma}

\begin{proof}
The elements of $\Int_*A$ are the finite unions of subsets of the
form $U=(b]_{\bullet}-X$, where $b\in A^-$ and $X$ is a finite union
of principal ideals of $A^-$. For
such a subset $U$, the properties (i) and (ii) above are obvious, with
$\Max U=\es$, if $X=(b]_{\bullet}$, $\Max U=\set{b}$, otherwise.
Furthermore, the set of all subsets $U$ of $A^-$ satisfying (i) and
(ii) above is obviously closed under finite union. Therefore, it
contains $\Int_*A$.
\end{proof}

\begin{definition}\label{D:StFn}
A map $\gx\colon A^-\to B$ is a \emph{step function}, if
the range of $\gx$ is finite, and the inverse
image $\gx^{-1}\set{b}$ belongs to $\Int_*A$, for all $b\in B^-$.
\end{definition}

If $\gx\colon A^-\to B$ is any map, a \emph{support} of $\gx$ is an
element $a$ of $A$ such that $\gx(x)>0$ implies $x\leq a$, for all
$x\in A^-$. By Lemma~\ref{L:MaxInt}, every step function has a support.

The following lemma establishes a useful compactness property of step
functions:

\begin{lemma}[Compactness of step functions]\label{L:CompStep}
 Let $\gh\colon A^-\to B$ be a step function, let $\E{D}$ be an upward
directed set of \emph{antitone} maps from $A^-$ to $B$. If for all $a\in
A^-$, there exists $\gx\in\E{D}$, with $\gh(a)\leq\gx(a)$, then there
exists an antitone map $\gx\in\E{D}$ such that $\gh\leq\gx$.
 \end{lemma}

\begin{proof}
 Put $U=\UUm{\Max\gh^{-1}\set{b}}{b\in\gh[A^-]-\set{0}}$. Since $\gh$ is
a step function, it follows from Lemma~\ref{L:MaxInt} that $U$ is finite.
For all $a\in U$, there exists, by assumption, an element $\gx_a$ of
$\E{D}$ such that $\gh(a)\leq\gx_a(a)$. Since $\E{D}$ is upward directed
and $U$ is finite, there exists $\gx\in\E{D}$ such that $\gx_a\leq\gx$,
for all $a\in U$. Now let $a\in A^-$; we prove that $\gh(a)\leq\gx(a)$.
This is trivial if $\gh(a)=0$. Now assume that $\gh(a)>0$. Since
$\gh^{-1}\set{\gh(a)}$ belongs to $\Int_*A$, there exists, by
Lemma~\ref{L:MaxInt}(i), $x\in\Max\gh^{-1}\set{\gh(a)}$ such that $a\leq
x$. Note that $x\in U$. Therefore,
 \begin{align*}
 \gh(a)&=\gh(x)&&(\text{since }x\in\gh^{-1}\set{\gh(a)})\\
 &\leq\gx_x(x)&&(\text{since }x\in U)\\
 &\leq\gx(x)&&(\text{by the definition of }\gx)\\
 &\leq\gx(a)&&
 (\text{since }a\leq x\text{ and }\gx\text{ is antitone}).
 \end{align*}
This holds for all $a\in A^-$, thus $\gh\leq\gx$.
\end{proof}

The following lemma provides us with a large supply of step functions:

\begin{lemma}\label{L:TensStep}
Let $n\in\go$, let $a_0$,\dots, $a_{n-1}\in A$, let $b_0$,\dots,
$b_{n-1}\in B$. Then the map $\gx\colon A^-\to B$ defined by
 \begin{equation}\label{Eq:TensStep}
 \gx(x)=\JJm{b_i}{i<n,\,x\leq a_i},
 \end{equation}
for all $x\in A^-$, is an antitone step function. Furthermore,
 \[
 \vv<x,\gx(x)>\in\JJm{a_i\otimes b_i}{i<n},
 \]
for all $x\in A^-$.
\end{lemma}

\begin{proof}
Everything in the statement of the lemma is obvious except for the fact
that $\gx$ is a step function. Put
 \[
 S(x)=\setm{i<n}{x\leq a_i},
 \]
for all $x\in A^-$. If $b\in B^-$, then, for all $x\in A^-$, $\gx(x)=b$
if{f} there exists a nonempty subset $I$ of $n$ such that
$
 \JJm{b_i}{i\in I}=b\text{ and }I=S(x).
$
Therefore, to prove that $\gx^{-1}\set{b}$ belongs to $\Int_*A$, it
suffices to prove that the set
 \[
 X_I=\setm{x\in A^-}{I=S(x)}
 \]
belongs to $\Int_*A$, for all nonempty $I\ci n$. But it is easy to verify
that
 \[
 X_I=\left(\MMm{a_i}{i\in I}\right]_{\bullet}-
 \UUm{(a_j]_{\bullet}}{j\in n-I},
 \]
which belongs to $\Int_*A$ since $I$ is nonempty.
\end{proof}

\subsection{The adjustment sequence of a step function}

The basic technical tool of transferability, see
Section~\ref{S:Transferability}, is the adjustment sequence: we adjust a
map to closer reflect the structure. Similar ideas come up in connection
with projective lattices (B. J\'onsson, first published in \cite{CJN77})
and bounded homomorphisms (R.~N. McKenzie \cite{rM70}); see also the
adjustment sequence of the $n$-modular identity in \cite{GrWe2}. The
following definition of the adjustment sequence takes into account that
the lattices may be \emph{infinite}. Note that we can only adjust maps
with finite range.

For all $a\in A^-$, we shall denote by $\CC(a)$ the set of all
nonempty, finite subsets $X$ of $A^-$ such that $a\leq\JJ X$.

\begin{definition}\label{D:OneStep}
 Let $\gx\colon A^-\to B$ be a map with finite range. The \emph{one-step
adjustment} of $\gx$ is  $\gx^{(1)}\colon A^-\to B$ defined by
 \begin{equation}\label{Eq:AdjSeq}
 \gx^{(1)}(x)=\JJ\left(\,\MM\gx[S]\bigm|S\in\CC(x)\,\right),
 \end{equation}
for all $x\in A^-$.
\end{definition}

Note that since the range of $\gx$ is finite, the right hand side of
the equation \eqref{Eq:AdjSeq} is well-defined.

\begin{remark}
Since $\set{x}\in\CC(x)$, for all $x\in A^-$, the inequality
$\gx\leq\gx^{(1)}$ always holds.
Let us assume, in addition, that $\gx$ is \emph{antitone}. Then the
expression \eqref{Eq:AdjSeq} for $\gx^{(1)}(x)$ takes on the following
form:
 \begin{equation}\label{Eq:AdjSeq2}
 \gx^{(1)}(x)=\gx(x)\jj\JJ\left(\,\MM\gx[S]\bigm|S\in\CC^*(x)\,\right),
 \end{equation}
where $\CC^*(x)$ denotes the set of all $S\in\CC(x)$ such that
$x\nin(S]$, see Section~\ref{S:LowBded}.

Let us further assume that $A$ is \emph{finite}. Then every
map $\gx\colon A^-\to B$ is a step function, and the one-step adjustment
of $\gx$ takes, on the join-irreducible elements of $A$, the following
form:
 \[
  \gx^{(1)}(p)=
 \gx(p)\jj\JJ\left(\,\MM\gx[I]\bigm| I\in\M(p)\,\right),
 \]
for all $p\in\J(A)$.
\end{remark}

\begin{lemma}\label{L:StCov}
Let $\gx\colon A^-\to B$ be a step function. For every subset
$X$ of $\gx[A^-]$, we define a subset
$\Cov(\gx;X)$ of $A^-$ as follows:
 \[
 \Cov(\gx;X)=\setm{x\in A^-}{\gx[S]=X,\text{ for some }S\in\CC(x)}
 \]
Then for all $a\in A$, $\Cov(\gx;X)\ii(a]_{\bullet}$ belongs to
$\Int_*A$.
\end{lemma}

\begin{proof}
Let $R=\gx[A^-]$ and define $U_b=\gx^{-1}\set{b}\ii(a]_{\bullet}$,
for all $b\in R$. Since $\gx$ is a step function, $\gx^{-1}\set{b}$
belongs to $\Int_*A$ for all $b\in R-\set{0}$, while
$\gx^{-1}\set{0}$ belongs to $\Int A$. Thus, $U_b$ belongs to
$\Int_*A$, for all $b\in R$.

Furthermore, $U=\UUm{\Max U_b}{b\in R}$ is a finite set, by
Lemma~\ref{L:MaxInt}(ii). Let $x\in\Cov(\gx;X)\ii(a]_{\bullet}$, that is,
$x\in A^-$, $x\leq a$, and there exists $S\in\CC(x)$ such that
$\gx[S]=X$. For all $s\in S$, $U_{\gx(s)}$ belongs to $\Int_*A$ and it
contains $s$ as an element. By Lemma~\ref{L:MaxInt}(i), there exists
$s^*\in\Max U_{\gx(s)}$ such that $s\leq s^*$. Put $T=\setm{s^*}{s\in
S}$. Then $T\in\CC(x)$, $\gx[T]=X$, and $T\ci U$. Hence, we have proved
the equality
 \[
 \Cov(\gx;X)\ii(a]_{\bullet}=\setm{x\in(a]_{\bullet}}
 {\gx[T]=X,\text{ for some }T\in\CC(x)\text{ with }T\ci U},
 \]
which can be written, by the definition of $\CC(x)$, as
 \[
 \Cov(\gx;X)\ii(a]_{\bullet}=
 \UU\left(\,\left(a\mm\JJ T\right]_{\bullet}\Bigm|
 T\ci U,\ \gx[T]=X\,\right).
 \]
Since $U$ is finite, $\Cov(\gx;X)\ii(a]_{\bullet}$ belongs to $\Int_*A$.
\end{proof}

\begin{lemma}\label{L:AdjHom}
Let $\gx\colon A^-\to B$ be a map with finite range. Then $\gx^{(1)}$
is antitone and has finite range. Furthermore, $\gx^{(1)}=\gx$
if{f} $\gx$ is a semilattice homomorphism from $\vv<A^-;\jj>$ to
$\vv<B;\mm>$.
\end{lemma}

\begin{proof}
Since $x\leq y$ implies that $\CC(y)\ci\CC(x)$, $\gx^{(1)}$ is obviously
antitone.
If $\gx$ is a semilattice homomorphism from $\vv<A^-;\jj>$ to
$\vv<B;\mm>$, then $\MM\gx[S]\leq\gx(x)$, for all $x\in A^-$ and all
$S\in\CC(X)$, thus $\gx^{(1)}(x)\leq\gx(x)$. Since $\gx\leq\gx^{(1)}$,
we obtain that $\gx=\gx^{(1)}$. Conversely, suppose that
$\gx=\gx^{(1)}$. Let $x$,
$y\in A^-$. Since $\set{x,y}\in\CC(x\jj y)$, the inequality
$\gx^{(1)}(x\jj y)\geq\gx(x)\mm\gx(y)$ holds. Since $\gx=\gx^{(1)}$ and
$\gx^{(1)}$ is antitone, we obtain that $\gx(x\jj y)=\gx(x)\mm\gx(y)$.
\end{proof}

\begin{lemma}\label{L:SuppOS}
Let $\gx\colon A^-\to B$ be a map with finite range. Then $\gx$ and
$\gx^{(1)}$ have the same supports.
\end{lemma}

\begin{proof}
Since $\gx\leq\gx^{(1)}$, every support of $\gx^{(1)}$ is a support of
$\gx$. Conversely, let $a$ be a support of $\gx$. Let $x\in A^-$ such
that
$x\nleq a$. For all $S\in\CC(x)$, there exists $s\in S$ such that
$s\nleq a$ (otherwise, $x\leq\JJ S\leq a$), thus $\MM\gx[S]=\gx(s)=0$;
whence $\gx^{(1)}(x)=0$. Hence $a$ is a support of $\gx^{(1)}$.
\end{proof}

We now state the main result of this section:

\begin{proposition}\label{P:StepAdj}
Let $\gx\colon A^-\to B$ be a step function. Then the one-step
adjustment $\gx^{(1)}$ of $\gx$ is an antitone step function.
\end{proposition}

\begin{proof}
The fact that $\gx^{(1)}$ is antitone with finite range has been
established in Lemma~\ref{L:AdjHom}.
So, to complete the proof that $\gx^{(1)}$ is a step function, it
suffices to prove that $(\gx^{(1)})^{-1}\set{b}$ belongs to $\Int_*A$,
for all
$b\in B^-$. Let $a_0$ be a support of $\gx$. Then, by
Lemma~\ref{L:SuppOS}, $a_0$ is also a support of $\gx^{(1)}$.
Put $R=\gx[A^-]$, and let $I$ be the set of
all nonempty subsets of $\Pow(R)-\set{\es}$.
Note that $I$ is finite. For all $\fa\in I$, define
 \begin{gather*}
 b_\fa=\JJ\left(\,\MM X\Bigm|X\in\fa\,\right),\\
 U_\fa=\setm{x\in(a_0]_{\bullet}}{\fa=\setm{\gx[S]}{S\in\CC(x)}}.
 \end{gather*}
We claim that the following equality holds:
 \begin{equation}\label{Eq:gx'-1b}
 (\gx^{(1)})^{-1}\set{b}=\UUm{U_\fa}{\fa\in I,\,b_\fa=b}.
 \end{equation}
Indeed, if $x$ belongs to the right hand side of \eqref{Eq:gx'-1b},
then there exists $\fa\in I$ such that $b_\fa=b$ and $x\in U_\fa$.
Hence $\fa=\setm{\gx[S]}{S\in\CC(x)}$, so that
 \begin{equation}\label{Eq:JJSXb}
 \gx^{(1)}(x)=\JJ\left(\,\MM\gx[S]\Bigm|S\in\CC(x)\,\right)
 =\JJ\left(\,\MM X\Bigm|X\in\fa\,\right)=b_\fa=b.
 \end{equation}
Conversely, let us assume that $\gx^{(1)}(x)=b$. Since $b>0$ and
$a_0$ is a support of $\gx^{(1)}$, it follows that $x\leq a_0$ and
$\fa=\setm{\gx[S]}{S\in\CC(x)}$ belongs to $I$, so that $x\in U_\fa$, by
definition. By an argument similar to \eqref{Eq:JJSXb}, we see that
$b_\fa=b$. This completes the proof of \eqref{Eq:gx'-1b}.

To complete the proof of Proposition~\ref{P:StepAdj}, it suffices
to prove that $U_\fa$ belongs to $\Int_*A$, for all $\fa\in I$. By
the definition of $U_\fa$, an element $x$ of $A^-$ belongs to $U_\fa$
if{f} the following three conditions are satisfied:
 \begin{enumerate}
 \item $x\leq a_0$;
 \item for all $X\in\fa$, there exists
$S\in\CC(x)$, satisfying $\gx[S]=X$;
 \item $\gx[S]\ne Y$ holds for all $Y\in\Pow(R)-\fa$ and $S\in\CC(x)$.
 \end{enumerate}
Therefore, we obtain that
 \[
 U_\fa = \IIm{(a_0]_{\bullet}\ii\Cov(\gx;X)}{X\in\fa}
 -\UUm{(a_0]_{\bullet}\ii\Cov(\gx;Y)}{Y\in\Pow(R)-\fa},
 \]
which belongs to $\Int_*A$, by Lemma~\ref{L:StCov}.
\end{proof}

By Proposition~\ref{P:StepAdj}, we can state the following definition:

 \begin{definition}\label{D:AdjSeq}
 Let $\gx\colon A^-\to B$ be a step function. The \emph{adjustment
sequence} of $\gx$ is the sequence $\vv<\gx^{(n)} \mid n\in\go>$
defined inductively by $\gx^{(0)}=\gx$, and
$\gx^{(n+1)}=(\gx^{(n)})^{(1)}$, for all $n\in\go$.
 \end{definition}

As an immediate consequence of Proposition~\ref{P:StepAdj}, we obtain
the following:

\begin{corollary}\label{C:AdjSeq}
Let $\gx\colon A^-\to B$ be a step function. Then the adjustment
sequence of $\gx$ is increasing, that is, $\gx^{(n)}\leq\gx^{(n+1)}$
for all $n$. Furthermore, if $n>0$, then $\gx^{(n)}$ is an
antitone step function.
\end{corollary}

\section{Capped tensor products and homomorphisms to
ideals}\label{S:ABSm}

In this section, we establish equivalent conditions under which a
tensor product of two lattices with zero is capped. These conditions
are stated in terms of adjustment sequences and choice functions (see
Section~\ref{S:homomorphisms}).

\begin{theorem}\label{T:CharABSm}
Let $A$ and $B$ be lattices with zero. Then the following statements are
equivalent:
\begin{enumerate}
\item $A\otimes B$ is a capped tensor product.

\item Let $\gf\colon\vv<A^-;\jj>\to\vv<\Id B;\ii>$ be a semilattice
homomorphism, let $\gx\colon A^-\to B$ be a step function. If $\gx$ is
a choice function for $\gf$, then there exists a step function
$\gh\colon A^-\to B$ such that the following conditions hold:
\begin{itemize}
\item[(ii${}_1$)] $\gx(a)\leq\gh(a)\in\gf(a)$, for all $a\in A^-$.

\item[(ii${}_2$)] The map $\gh$ is a semilattice homomorphism from
$\vv<A^-;\jj>$ to $\vv<B;\mm>$.
\end{itemize}

\item The adjustment sequence of any step function
$\gx\colon A^-\to B$ is eventually constant.
\end{enumerate}
\end{theorem}

\begin{proof}
Throughout this proof, we shall denote by
$\ge\colon A\ootimest B\to A\ootimes B$ the
canonical isomorphism (see Section~\ref{S:antitone}).

(i) \emph{implies} (ii). We claim that
 \[
  H=\JJm{a\otimes\gx(a)}{a\in A^-}\in A\ootimes B
 \]
 belongs to $A\otimes B$. Indeed, it follows from
Lemma~\ref{L:MaxInt}(i) that the equality
 \[
   \JJm{a\otimes b}{a\in\gx^{-1}\set{b}}=
   \JJm{a\otimes b}{a\in\Max\gx^{-1}\set{b}}
 \]
 holds, for all $b\in\gx[A^-]-\set{0}$. Thus, the equality
 \begin{equation}\label{Eq:Hfin}
    H=\JJm{a\otimes\gx(a)}
    {b\in\gx[A^-]-\set{0},\,a\in\Max\gx^{-1}\set{b}}\in A\otimes B
 \end{equation}
holds, and the right hand side of \eqref{Eq:Hfin} is a finite join, by
Lemma~\ref{L:MaxInt}(ii). Hence $H$ belongs to
$A\otimes B$.

Since $A\otimes B$ is a capped tensor product, $H$ is a finite union of
pure tensors. Thus all the values of the map $\ga=\ge^{-1}(H)$ are
principal ideals of $B$, say, $\ga(a)=(\gh(a)]$, for all $a\in A^-$. Since
$\ga\in A\ootimest B$, the map $\gh$ is a homomorphism from
$\vv<A^-;\jj>$ to $\vv<B;\mm>$. Furthermore, $\vv<a,\gx(a)>\in H$ for
all $a\in A^-$, thus $\gx\leq\gh$. Put $K=\ge(\gf)$. For all
$a\in A^-$, the pair $\vv<a,\gx(a)>$ belongs to $K$; thus $H\ci K$. It
follows that $\ge^{-1}(H)\leq\ge^{-1}(K)$, that is, $\gh(a)\in\gf(a)$,
for all $a\in A^-$. Thus $\gh$ is a choice function for $\gf$.

(ii) \emph{implies} (iii). Let $\gx\colon A^-\to B$ be any step
function. We associate with $\gx$ its adjustment sequence,
$\vv<\gx^{(n)} \mid n\in\go>$. For all $a\in A^-$, we define $\gf(a)$ as
the ideal of $B$ generated by the set $\setm{\gx^{(n)}(a)}{n\in\go}$. By
Corollary~\ref{C:AdjSeq}, all the maps $\gx^{(n)}$ are
antitone, for $n>0$, thus $x\leq y$ implies $\gf(x)\ce\gf(y)$, for $x$,
$y\in A$. Furthermore, if $x$, $y\in A^-$, then, since
$\set{x,y}\in\CC(x\jj y)$, the inequality
$\gx^{(n+1)}(x\jj y)\geq \gx^{(n)}(x)\mm\gx^{(n)}(y)$
holds, for all $n$. Therefore, $\gf(x \jj y) \ce \gf(x)
\ii \gf(y)$; since $\gf$ is antitone, $\gf(x \jj y) = \gf(x) \ii \gf(y)$.
Thus $\gf$ is a semilattice homomorphism from $\vv<A^-;\jj>$ to $\vv<\Id
B;\ii>$. Since $\gx$ is a choice function for $\gf$, there exists, by
assumption, a step function $\gh\colon A^-\to B$ such that $\gx\leq\gh$
and $\gh$ is a choice function for $\gf$.

For all $a\in A^-$, by the definition of $\gf$, there exists $n>0$
such that $\gh(a)\leq\gx^{(n)}(a)$. For $n>0$, since the $\gx^{(n)}$
are antitone and since $\gh$ is a step function, it follows from
Lemma~\ref{L:CompStep}, applied to $\E{D}=\setm{\gx^{(n)}}{n>0}$,
that there exists $m>0$ such that $\gh\leq\gx^{(m)}$. On the other
hand, $\gx\leq\gh$ and $\gh$ is a homomorphism from $\vv<A^-;\jj>$ to
$\vv<B;\mm>$, thus, by Lemma~\ref{L:AdjHom},
$\gx^{(n)}\leq\gh^{(n)}=\gh$, for all $n\in\go$. It follows that
$\gh=\gx^{(m)}=\gx^{(m+1)}$.

(iii) \emph{implies} (i).
Let $H \in A \otimes B$. We prove that $H$ is
a finite union of pure tensors. Write
 \[
 H=\JJm{a_i\otimes b_i}{i<n},
 \]
where $n\in\go$ and $\vv<a_i,b_i>\in A\times B$, for all $i<n$.
Consider the function $\gx\colon A^-\to B$ given by the formula
\eqref{Eq:TensStep}. By Lemma~\ref{L:TensStep}, $\gx$ is an antitone
step function. By~assumption, there exists $m\in\go$ such that
$\gx^{(m)}=\gx^{(m+1)}$. Put $\gh=\gx^{(m)}$. By~Lemma~\ref{L:AdjHom},
$\gh$ is a homomorphism from $\vv<A^-;\jj>$ to $\vv<B;\mm>$. In
particular, the set
 \[
 K=\UUm{a\otimes\gh(a)}{a\in A^-}
 \]
is a bi-ideal of $A\times B$. Thus, since $\gx\leq\gh$
and $\gx(a_i)\geq b_i$ for all $i<n$, $K$ contains $H$.

Put $\gf=\ge^{-1}(H)$. By Lemma~\ref{L:TensStep},
$\vv<x,\gx(x)>\in H$ for all $x\in A^-$, thus $\gx$ is a choice
function for $\gf$. Since $\gf\in A\ootimest B$, all the $\gx^{(n)}$,
$n\in\go$, are choice functions for $\gf$. In particular, $\gh$ is a
choice function for $\gf$. This means that $\vv<a,\gh(a)>\in H$, for
all $a\in A^-$, that is, $K$ is contained in $H$. Finally, $K=H$.

Since $\gh$ is both a step function and a homomorphism from
$\vv<A^-;\jj>$ to $\vv<B;\mm>$, the set $\gh^{-1}\set{b}$ has, by
Lemma~\ref{L:MaxInt}, a largest element, say, $a_b$, for all
$b\in\gh[A^-]-\set{0}$. It follows that
 \[
 H=K=\UUm{a_b\otimes b}{b\in\gh[A^-]-\set{0}},
 \]
a finite union of pure tensors.
\end{proof}

\section{amenability and \TJ}\label{S:TJimplies}
In this section, we characterize amenable lattices. We start with
a sufficient condition.

\begin{proposition}\label{P:TJimpSm}
Let $A$ be a finite lattice. If
$A$ satisfies \TJ, then $A$ is amenable.
\end{proposition}

\begin{proof}
 Let $A$ be a finite lattice satisfying \TJ. Since $A$ satisfies \TJ,
$\J(A)=\set{p_1,\ldots,p_n}$ so that $p_i\mathrel{D}p_j$ implies that
$i>j$, for all $i$, $j$ in $\set{1,\ldots,n}$. We~verify Condition (iii)
of Theorem~\ref{T:CharABSm} for $A$. So let $B$ be any lattice with zero,
and let $\gx\colon A^-\to B$ be any map. We prove that the adjustment
sequence of $\gx$ is eventually constant. By replacing $\gx$ by
$\gx^{(1)}$, we may assume, without loss of generality, that $\gx$ is
antitone.

\begin{claim}
Let $0<j<i\leq n+1$. Then $\gx^{(i)}(p_j)=\gx^{(j)}(p_j)$.
\end{claim}

\begin{proof}[Proof of Claim]
We prove the Claim by induction on $i$. The Claim is
vacuously true for $i=0$. Let $0<j<i+1\leq n+1$, and let us
assume that the induction hypothesis holds for $i$.

By definition,
 \begin{equation}\label{Eq:xi+1i}
 \gx^{(i+1)}(p_j) = \gx^{(i)}(p_j)\jj
 \JJ\left(\, \MM\gx^{(i)}[I]\mid I\in\M(p_j)\,\right).
 \end{equation}

If $j\leq 2$, then $\M(p_j)=\es$ (because if $\vv<p,I>$ is a
minimal pair, then $I$ has at least two elements), thus, by the
induction hypothesis,
 \[
 \gx^{(i+1)}(p_j) = \gx^{(i)}(p_j) = \gx^{(j)}(p_j).
 \]

If $j>2$, then let $I\in\M(p_j)$ and let $k$ satisfy $p_k\in I$. In
particular, $k<j$, thus $k<i$, and so by the induction hypothesis,
 \[
 \gx^{(i)}(p_k)=\gx^{(i-1)}(p_k)=\gx^{(k)}(p_k).
 \]
Therefore,
$\MM\gx^{(i)}[I]=\MM\gx^{(i-1)}[I]$. Thus, applying
\eqref{Eq:xi+1i} to $i-1$, we obtain that
$\MM\gx^{(i)}[I]\leq\gx^{(i)}(p_j)$. Using \eqref{Eq:xi+1i} and the
induction hypothesis,
 \begin{equation*}
 \gx^{(i+1)}(p_j) = \gx^{(i)}(p_j) = \gx^{(j)}(p_j),
 \end{equation*}
completing the proof of the Claim.
\end{proof}

By applying the Claim for $i=n+1$, we obtain that $\gx^{(n)}$ and
$\gx^{(n+1)}$ agree on $\J(A)$, thus $\gx^{(n)}=\gx^{(n+1)}$.
We have verified Theorem~\ref{T:CharABSm}(iii) for $A$,
so $A$ is amenable.
\end{proof}

Now we are ready for the characterization of amenable lattices:

\begin{theorem}\label{T:JoinAm}
Let $A$ be a lattice with zero. Then the following conditions are
equivalent:
 \begin{enumerate}
 \item $A$ is amenable.
 \item $A$ is locally finite and every finite sublattice of $A$
satisfies \TJ.
 \end{enumerate}
 \end{theorem}

 \begin{proof}
 (ii) \emph{implies} (i). Let $L$ be a lattice with zero. By
Proposition~\ref{P:TJimpSm}, $U\otimes L$ is a capped tensor product, for
every finite $\set{0}$-sublattice $U$ of $A$. By
Proposition~\ref{P:SmLim}, $A\otimes L$ is a capped tensor product.

(i) \emph{implies} (ii). It suffices to prove that for every finitely
generated amenable lattice $A$ with zero, $A$ is finite and it
satisfies \TJ. Let $X$ be a finite generating subset of $A$, and let
$h\colon\FL(X)\twoheadrightarrow A$ be the canonical surjective
homomorphism. We~shall consider the lower limit table
$\vv<\gb_n \mid n\in\go>$ associated with $h$ and $X$. The maps $\gb_n$
can be computed using the formulas \eqref{Eq:beta0} and \eqref{Eq:betan}.

Now we shall use the fact that $A\otimes\FL(X)$ is a capped tensor
product. We consider the step function $\gx\colon A^-\to\FL(X)$
defined by \eqref{Eq:TensStep}, see Lemma~\ref{L:TensStep}. That is,
 \[
   \gx(a)=\JJm{x\in X}{a\leq x},
 \]
 for all $a\in A^-$.  Let $w \mapsto w^\dd$ denote the dualization map on
$\FL(X)$. We obtain, in particular, that $\gx(a)=\gb_0(a)^\dd$.

 Now we consider the adjustment sequence, $\vv<\gx^{(n)} \mid n\in\go>$,
of $\gx$. Let $n\in\go$, and let us assume that we have proved that
$\gx^{(n)}(a)=\gb_n(a)^\dd$, for all $a\in A^-$. We~compute
$\gx^{(n+1)}(a)$, for $a\in A^-$:
 \begin{align*}
 \gx^{(n+1)}(a)&=\gx^{(n)}(a)\jj\JJ\left(\,
 \MM\gx^{(n)}[S]\mid S\in\CC(a)\,\right)\\
 &=\gb_n(a)^\dd\jj\JJ\left(\,
 \MM\gb_n[S]^\dd\mid S\in\CC(a),\,\JJ S\leq h(1)\,\right)\\
\intertext{(the condition $\JJ S\leq h(1)$ is satisfied, because
$h(1)=1$)}
 &=\gb_{n+1}(a)^\dd,
 \end{align*}
the last step by \eqref{Eq:betan}. It follows that $\gx^{(n)}(a) =
\gb_n(a)^\dd$, for all $n\in\go$ and all $a\in A^-$. However, since
$A\otimes\FL(X)$ is a capped tensor product, there exists, by
Theorem~\ref{T:CharABSm}, $n\in\go$ such that $\gx^{(n)}=\gx^{(n+1)}$. It
follows that $\gb_n=\gb_{n+1}$. By Proposition~\ref{P:LoBdn}, $A$ is
finite and $h$ is a lower bounded homomorphism. By
Proposition~\ref{P:LBfin}, $A$ satisfies \TJ.
 \end{proof}

The smallest finite lattice that does not satisfy \TJ\ is the diamond
$M_3$. Since $M_3$ is a quotient of $\FL(3)$, one obtains, by
Corollary 3.7 of \cite{GrWe}, the following:

 \begin{corollary}\label{C:TensNotLatt}
Neither $M_3 \otimes\FL(3)$
nor $\FL(3) \otimes \FL(3)$ is a lattice.
 \end{corollary}

\begin{corollary}\label{C:SmNonVar}
There exists a countable sequence
$\vv<S_n\mid n\in\go>$ of finite amenable lattices such that the
product $\prod_{n\in\go}S_n$ is not amenable.
\end{corollary}

\begin{proof}
Denote by $x_0$, $x_1$, and $x_2$ the generators of the
free lattice $\FL(3)$ on three generators. The subset $S_n = S_n(3)$
(see Section~\ref{S:set}) of $\FL(3)$ is a finite \jz-subsemilattice
of $\FL(3)$, thus it is a lattice. Furthermore, by the end of the
proof of Lemma 2.77 in \cite{FJNa95}, all $S_n$ are lower bounded.
Therefore, by Proposition~\ref{P:LBfin}, all $S_n$ satisfy \TJ.
By Proposition~\ref{P:TJimpSm}, all $S_n$ are amenable.

On the other hand, the diagonal map embeds $\FL(3)$ into the reduced
product $L = \prod_{\E F}\vv<S_n\mid n \in\go>$, where $\E F$ denotes
the Fr\'echet filter on $\go$. Suppose that $S = \prod_{n \in\go}S_n$
is amenable. Since $L$ is a quotient of $S$ and $\FL(3)$ embeds
into~$L$, it follows from Theorem~\ref{T:SmPres} that $\FL(3)$ is
also amenable, a contradiction by Corollary~\ref{C:TensNotLatt}.
\end{proof}

In fact, the proof above shows that $M_3\otimes S$ is not a lattice.

\begin{example}\label{E:LocFin}
 The proof of Corollary~\ref{C:SmNonVar} gives immediately a locally
finite \emph{amenable} lattice $S$ with zero such that $M_3$ belongs to
the variety generated by $S$ (in fact, $S$ generates the variety $\mbf
L$ of all lattices): consider the semilattice direct sum $\bigoplus_{n
\in\go}S_n$, where the $S_n$ are the finite lattices in the proof of
Corollary~\ref{C:SmNonVar}. This shows that the hypothesis of
Corollary~\ref{C:SmFinVar} that $B$ is \emph{finite} cannot be
weakened to $B$ being locally finite.
 \end{example}

In particular, the class of amenable lattices is not a variety.

\section{Discussion}\label{S:problems}

The most central open question is stated first:

\begin{problem}\label{Q:EquiSm}
 Let $A$ be a lattice with zero. If, for every lattice $L$ with zero,
$A\otimes L$ is a lattice, is $A$ amenable?
 \end{problem}

By Theorem~\ref{T:SmFin}, the answer to Problem~\ref{Q:EquiSm} is
positive for  locally finite lattice $A$. We believe that in the
general case the answer is in the negative.

 \begin{problem}\label{Pb:SimAme}
 Does there exist a nontrivial, simple, amenable lattice with zero?
 \end{problem}

As we will show in Corollary~\ref{C:SimpSDj}, there is no
nontrivial simple amenable lattice with a largest element.

We recall that a lattice $L$ is \emph{join-semidistributive}, if it
satisfies the following condition:
 \begin{equation}
 x\jj z=y\jj z\text{\q implies that\q }x\jj z=(x\mm y)\jj z,
 \q\text{for all }x,\,y,\,z\in L.\tag*{\SD}
 \end{equation}

\begin{proposition}\label{P:SimpSDj}
Let $S$ be a simple lattice with at least three elements. If $S$
satisfies \SD, then $S$ does not have a largest element.
\end{proposition}

\begin{proof}
If $S$ has a largest element, then it has a maximal ideal, say, $I$.
Then $I$ is a \emph{prime ideal} of $S$. Indeed, if $x$, $y\nin I$, then,
by the maximality of $I$, there exists $u\in I$ such that $x\jj u=y\jj u=1$.
By \SD, $(x\mm y)\jj u=1$, thus $x\mm y\nin I$ (otherwise,
$1\in I$, a contradiction).

So $I$ defines a lattice homomorphism
from $S$ to the two-element chain. Since $S$ has at least three elements,
the kernel of this homomorphism is a non-trivial congruence of $S$, which
contradicts the simplicity of~$S$.
\end{proof}

It is known that every finite lattice satisfying \TJ\ satisfies \SD,
see, for example, Theorem~2.20 in \cite{FJNa95}. Since \SD\ is
preserved under direct limits, we obtain the following corollary:

\begin{corollary}\label{C:SimpSDj}
Let $S$ be a simple amenable lattice with zero. If $S$ has at least
three elements, then $S$ does not have a largest element.
\end{corollary}

To formulate the next problem, let us introduce an additional
terminology.

\begin{definition}\label{D:CJA}
Let $A$ be a lattice with zero and let $\mbf{C}$ be
a class of lattices. Then $A$ is \emph{$\mbf{C}$-amenable}, if $A
\otimes L$ is a lattice, for every lattice with zero $L$ in $\mbf{C}$.
\end{definition}

\begin{problem}\label{Q:modcap}
Let $\mbf{V}$ be a variety of
lattices. Is the class of finite $\mbf{V}$-amenable lattices
\emph{decidable}?
\end{problem}

For example, if $\mbf{M}$ is the variety of all \emph{modular}
lattices, then $M_3$ is $\mbf{M}$-amenable. It would be desirable to
obtain a combinatorial characterization of $\mbf{M}$-amenable
lattices, as we did in this paper for $\mbf{L}$-amenable lattices,
where $\mbf{L}$ is the variety of all lattices. On the other extreme,
if $\mbf{D}$ is the variety of all distributive lattices, every
lattice with zero is $\mbf{D}$-amenable.

\begin{problem}
 Let $\mbf{A}$ denote the variety of all \emph{Arguesian} lattices. For a
finite lattice~$A$, prove that $A$ is $\mbf{M}$-amenable if{f} it is
$\mbf{A}$-amenable.
 \end{problem}

By using a lattice constructed in \cite{DaHW72}, we can prove that
$M_4$ is not $\mbf{M}$-amenable, see \cite{GrWe2} for details.
Furthermore, the corresponding counterexample is a lattice of
subspaces of a vector space, thus it is Arguesian. So one may expect,
for every non $\mbf{M}$-amenable lattice $A$, the existence of a
lattice $L$ of subspaces of some vector space such that $A \otimes L$
is not a lattice.

We have proved in Theorem~\ref{T:JoinAm} that every amenable
lattice is locally finite. This result cannot be relativized to
arbitrary varieties, as, for example, any lattice with zero is
$\mbf{D}$-amenable.

\begin{problem}
Is every $\mbf{M}$-amenable lattice locally finite?
\end{problem}

\section*{Acknowledgment} This work was partially completed while the
second author was visiting the University of Manitoba. The excellent
conditions provided by the Mathematics Department, and, in particular,
a quite lively seminar, were greatly appreciated.

The authors wish to thank the referee for some very constructive
suggestions.

\end{document}